\documentclass[11pt,twoside]{article}

\usepackage{color}
\usepackage{amsfonts}
\usepackage{amsmath}
\usepackage{amssymb}

\newcommand{\bean}{\begin{eqnarray*}}
\newcommand{\eean}{\end{eqnarray*}}
\newcommand{\benu}{\begin{enumerate}}
\newcommand{\eenu}{\end{enumerate}}
\newcommand{\eea}{\end{eqnarray}}
\newcommand{\bea}{\begin{eqnarray}}

\newtheorem{Theorem}{Theorem}

\newtheorem{Lemma}{Lemma}

\newcommand{\be}{\begin{equation}}
\newcommand{\ee}{\end{equation}}

\newcommand{\N}{{\mathbb N}}
\newcommand{\Z}{{\mathbb Z}}
\newcommand{\R}{{\mathbb R}}

\newcommand{\SSS}{{\mathbb S}}

\newcommand{\la}{\lambda}

\newcommand{\e}{\epsilon}
\newcommand{\1}{\mathbf{1}}

\def\g{\gamma}
\def\d{\delta}
\def\l{\lambda}
\def\o{\circ}

\newcommand{\ben}{\begin{enumerate}}
\newcommand{\een}{\end{enumerate}}
\newcommand{\bit}{\begin{itemize}}
\newcommand{\eit}{\end{itemize}}
\newcommand{\edoc}{\end{document}}

\newcommand{\bdefi}{\begin{definition}}
\newcommand{\btheo}{\begin{theorem}}
\newcommand{\bprop}{\begin{proposition}}
\newcommand{\brema}{\begin{remark}}
\newcommand{\bcoro}{\begin{corollary}}
\newcommand{\blemm}{\begin{lemma}}
\newcommand{\bexam}{\begin{example}}

\newcommand{\edefi}{\end{definition}}
\newcommand{\etheo}{\end{theorem}}
\newcommand{\eprop}{\end{proposition}}
\newcommand{\erema}{\end{remark}}
\newcommand{\ecoro}{\end{corollary}}
\newcommand{\elemm}{\end{lemma}}
\newcommand{\eexam}{\end{example}}

\newcommand{\V}{\noindent}
\newcommand{\ci}{\circ}

\newcommand{\Gl}{{\rm Gl}}
\newcommand{\Hol}{{\rm Hol}}
\newcommand{\Met}{{\rm Met}}
\newcommand{\diam}{{\rm diam}}

\newcommand{\cl}{{\rm cl}}

\newcommand{\Riem}{{\rm Riem}}

\def\qed{\ensuremath{\quad\Box\quad}}

\title{Restricted holonomy is lower semicontinuous}

\author{Olaf M\"uller\footnote{Humboldt-Universit\"at zu Berlin, Institut f\"ur Mathematik, Unter den Linden 6, 10099 Berlin. Email: mullerol@math.hu-berlin.de}}

\date{\today}

\begin{document}

\maketitle

\begin{abstract}
\V In this article, we examine continuity properties of the maps $\Hol$ and $\Hol^0$ assigning, on a fixed manifold $M$, to a metric on $M$ its holonomy class resp. restricted holonomy class (conjugacy class of the connected component of the holonomy representation). Among related results, we show that $\Hol^0$ is lower semicontinuous w.r.t. $C^1$ topology on the space of $C^2$ metrics.
\end{abstract}  

\section{Introduction and statement of the main result}

\V A vital feature of control theory and its innumerous technical applications (like parking a car or steering a satellite via 2 or 3 nozzles) is the effect of holonomy. Let $M$ be an oriented $n$-dimensional manifold. If we consider holonomy of the Levi-Civita connection as a map $\Hol$ on $\Met(M)$ to the set of conjugacy classes of subgroups of $SO(n)$, then an important question is whether $\Hol$ --- or at least its connected component $\Hol^0$ --- has some continuity properties. This article will show a partially affirmative answer. Let $M$ be an $n$-dimensional connected manifold and let $\Met (M)$ be the space of Riemannian metrics on $M$ of regularity\footnote{About the choice of this regularity, see the last paragraph of the article.} $C^2$, equipped with the usual $C^2 $ (metrizable) topology (i.e., a sequence of metrics converges iff every differential of order less or equal $2$ converges uniformly). Let $S(G)$ resp. $SC(G)$ resp. $SCC(G)$ be the set of conjugacy classes of subgroups resp. of subgroups with finitely many connected compact components resp. of compact connected subgroups of a Lie group $G$, and $P_n:= S(SO(n))$ and $PC_n:= SCC(SO(n)$). It is well-known that the maps $\Hol: \Met (M) \rightarrow S(SO(n))$, defined by $ \Hol(g)= [ \Hol_p (g) ] $, and $\Hol^0: \Met (M) \rightarrow PC_n $, with $ \Hol^0 (g) = [\Hol^0_p (g)] \forall g \in \Met (M)$ (where we use any $g$-orthormal frame to identify $T_pM$ with $\R^n$) do not depend on the choice of $p \in M$ (see e.g. \cite{KN}, Th. IV.5.5 for $\Hol^0$ taking values in $PC_n$).
 
\V On $SC(G)$ there is a (partial) order $\leq $ by inclusion of representatives via

\bea
\label{po-def}
&\forall a,b \in SC(G)&: \big( a \leq b : \Leftrightarrow \exists x \in a \exists y \in b: x \subset y \big) , {\rm \ more \ explicitly,} \nonumber \\
&[A] \leq [B]& : \Leftrightarrow \exists g \in G: gAg^{-1} \subset B,  
 \eea

\V cf. Theorem \ref{PartialOrder} (which makes $SCC(G)$ a poset by restriction, too).

\V A map $f$ from a topological space $X$ to a poset $Y$ is called {\bf lower (resp. upper) semicontinuous (l.s.c. resp. u.s.c.) at $x  \in X$} iff for all $y \in Y $ with $y < f(x) $ (resp. $y > f(x)$) there is a neighborhood $U$ of $x$ with $ f(u) \geq y \forall u \in U$ (resp. $f(u) \leq y \forall u \in U$). In Theorem \ref{semicontinuities} we will see that if $Y$ is a power set with the inclusion or $PC_n$ as above, where the map $\Hol^0$ takes values, a map $f: X \rightarrow Y  $ is lower semi-continuous at $x \in X$ if and only if there is a neighborhood $U$ of $x$ such that $f(x) \subset f(u)$ for all $u \in U$. The main results of this article is:

\begin{Theorem}
\label{main}
Let $M$ be a manifold and let $\pi: E \rightarrow M$ be a $G$-principal bundle. Denote with ${\rm Con} (\pi)$ the space of $G$-principal connections on $\pi$ of regularity $C^1$. 
\begin{enumerate}
	\item Let $ \nabla \in {\rm Con} (\pi)$. If there is a $C^0$-neighborhood $U \subset {\rm Con} (\pi)$ of $\nabla$ and $D \in \N$ s.t. the intrinsic diameter of ${\rm Hol}^0 (h)$ is $\leq D$, then ${\rm Hol}^0$ is l.s.c. at $ \nabla$. 
	\item $\Hol^0_p: {\rm Met} (M) \rightarrow PC_n$ is l.s.c. w.r.t. the $C^1$ metric on ${\rm Met} (M)$ at any $C^2$ metric $g$ on $M$.
\end{enumerate} 
\end{Theorem}

\V Wilking \cite{Wilking} found five-dimensional examples of compact Riemannnian manifolds of holonomies with countably infinitely many connected omponents. Moreover, $\Hol^0 (\nabla)$ is not compact for all linear connections $\nabla$: Hano and Ozuki \cite{jHhO} gave an example of a torsion-free connection on $\R^6$ whose holonomy is not closed in $SO(n)$, and there is a straightforward modification of that example that yields a metric connection on a rank-$4$ vector bundle over $\R^2$ whose holonomy is not closed in $SO(4)$. Moreover, the holonomies so obtained are subgroups of the form $G_a:= \{ R_\theta  \oplus R_{a \theta} | \theta \in \R \} \subset SO(4)$ where $ R_\theta $ is the $2$-dimensional rotation by the angle $\theta$, the direct sum denotes a block matrix in a fixed basis and we consider the limit $a \rightarrow 0$. The calculation of the characteristic polynomials shows that $G_a $ is conjugate to $G_b$ if and only if $a = b$. This shows that there is no statement similar to our main result for holonomies of metric connections on general vector bundles, nor for torsion-free connections on the tangent bundle. On the other hand, if the total space of the universal covering map $\pi$ of a compact manifold $M$ is spin, then for all $g \in \Met(M)$ such that $\pi^* g$ carries a non-zero parallel spinor, the holonomy class in the connected component of $g$ in $\Met (M)$ is even constant, cf. \cite{AKWW}.

\V The proof of the main result contains interesting side results, like a quantified and a converse Montgomery-Zippin Theorem. 

\bigskip

\V The author would like to thank Bernd Ammann for helpful discussions, and an anonymous referee for helpful comments.

\section{Proofs}

\V {\bf Conventions:} For a topological space $T$ and $x \in T$ let $N(x)$ be the set of neighborhoods of $x$, and let $C(T)$ be the set of closed subsets of $T$. All manifolds in this article are supposed to be finite-dimensional. Fix an oriented $n$-dimensional manifold $M$ throughout the article. The connected component of the identity in a topological group $G$ is called $G^0$. The Lie algebra of a Lie group $G$ is denoted by $LA(G)$, and, for $0 \leq k \leq n = \dim (V)$, the oriented $k$-Grassmannian of a vector space $V$ with a scalar product $\langle \cdot, \cdot \rangle $ by ${\rm Gr}_k(V)$, topologized by the metric $d_P(A,B) := \vert P_A - P_B \vert_{{\rm op}}$ where $A,B$ are two $k$-dimensional subspaces, $P_K$ is the orthonormal projection onto $K$ for any subspace $K$ of $V$, and $\vert \cdot \vert_{{\rm op}}$ is the operator norm of a linear map. 

\bigskip

\V We can equip the space $C_n $ of closed subgroups of $SO(n)$ with the quotient space topology under $q: C (SO(n)) \rightarrow C_n$ where, for a group $K$, the space $C(K)$ resp. $FC(K)$ is the set of closed subgroups of $K$ resp. of closed subgroups of $K$ with finitely many connected components, both equipped with the Hausdorff distance (defined from a bi-invariant metric on $SO(n)$) and its induced topology (and the identity component $\Hol^0_p (g)$ of $\Hol^0_p (g)$ {\em is} a closed subgroup of $SO(n)$, thus a compact Lie group).   

\V The small side question arises if we could also order the conjugacy classes of groups by requiring not only injective maps but injective group homomorphisms. This is indeed the case, as we will see now.

\V Let us first revise an example due to Tsit Yuen Lam showing that even in the category of subgroups of $\Gl(2, \R)$ up to conjugacy there is no Cantor-Bernstein theorem: Consider the usual linear representation $\Phi $ of the affine group in one dimension given by $\Phi: Y:= (\R \setminus \{ 0 \}) \times \R \rightarrow {\rm Gl}(2,\R) $ defined by $(x,y) \mapsto \left( \begin{array}{rr} x & y\\ 0&1 \end{array}\right)$. Now $G := \Phi (Y) $ is a subgroup of $\Gl (2, \R)$, and $H_z:= \Phi (\{ 1\} \times z \cdot \Z)$ is, for each $z \in \Z$, a subgroup of $G$ and of $\Gl(2, \R)$. For $g:= \Phi ((2,0))$ we get $g H_1 g^{-1} = H_2 \subsetneq H_1$.

\V In contrast to the absence of a Cantor-Bernstein Theorem in the entire category of groups, there is even a {\em strong} Cantor-Bernstein Theorem for compact Lie groups\footnote{"strong" means that if there are injective morphisms $f: X \rightarrow Y$ and $g : Y \rightarrow X$ then $f$ and $g$ {\em themselves} are isomorphisms.}:

\begin{Theorem}
	\label{PartialOrder}
	Let $K,H$ be compact Lie groups, let $f: K \rightarrow H$ and $g: H \rightarrow K$ be injective Lie group homomorphisms. Then $f$ and $g$ are Lie group isomorphisms. 
	
	\V Furthermore, let $G$ be a compact Lie group and let $K , H $ be closed subgroups of $G$ with finitely many connected components such that there are $g,h \in G$ with $g K g^{-1} \subset H$ and $h H h^{-1} \subset K$. Then $g K g^{-1} = H $. 
\end{Theorem}

\V {\bf Proof.} $d_1 f : T_1 K \rightarrow T_1 H $ and $d_1 g: T_1 H \rightarrow T_1 K$ are linear. Moreover, the fact that for the local diffeomorphism $\exp$ we have $\exp \ci d_1 f = f \o \exp$ implies that $d_{1} f$ is injective: Assume $d_{1} f (v) = 0 $ for some $v \in T_{1} K $, then choose $U \subset T_{1} K $ such that $\exp \vert_U $ is a diffeomorphism onto its image. There is $\lambda  \in (0; \infty)$ with $w:= \l v \in U$ and still $d_{1} f (w) = 0$. Now we have $1 = \exp (d_{1} f(w)) = f (\exp (w))$, so $\exp (w) =1$ and thus $w= 0$ and $v = 0$.

\V By the Cantor-Bernstein Theorem in the category of vector spaces, both tangent spaces have the same dimension and both linear maps are isomorphisms. Moreover, the connected component $H_0$ of $1$ in $H$ is contained in $f(K)$ by surjectivity of $\exp$ onto $H_0$ and surjectivity of $d_{\1} f$. With an analogous argument, $f$ is surjective onto every connected component it meets. Therefore $\pi_0 f : \pi_0 K \rightarrow \pi_0 H $ is injective. As both $\pi_0 K $ and $\pi_0 H$ are finite, they contain the same number of elements, and $\pi_0 f  $ is bijective. 

\V The second statement follows from the first by application to the Lie group isomorphism defined by conjugation with elements $g$ and $h$. \hfill \qed

\bigskip

Now to the question of semi-continuity. On the space $C(X)$ of closed subsets of a topological space $X$ we define a relation $\leq_C$ by $A \leq_C B :  \Leftrightarrow  A=B \lor A \subset {\rm int} (B)$, which is easily seen to be an order relation.

\begin{Theorem}
	\label{semicontinuities}
	Let always $f: T \rightarrow Y$ be a map from a topological space $T \ni x$. In Items 2, 3 recall that the cardinality of a conjugacy class is well-defined.
	\begin{enumerate} 
		\item Let $(Y, \leq )= (P(S), \subset) $ for a set $S$. If $\# f(x) \geq 2 $ resp. $ \# ( S \setminus f(x) ) \geq 2 $, then $f$ is lower resp. upper semi-continuous at $x$ if and only if there is $U \in N(x)$ with $f(u) \geq f(x) \forall u \in U$ resp. $f(u) \leq f(x) \forall u \in U$.
		\item For a Lie group $G$, let $Y = S(G) $ be the set of conjugacy classes of compact subgroups of $G$ with the order of Eq. \ref{po-def} (cf. Th. \ref{PartialOrder}). If $\# f(x) \geq 2$ resp. $\# (S \setminus f(x)) \geq 2$, then $f$ is lower resp. upper semi-continuous at $x$ if and only if there is $U \in N(x)$ with $f(u) \geq f(x) \forall u \in U$ resp. $f(u) \leq f(x) \forall u \in U$.
		\item Let $(M,d)$ be a Heine-Borel metric space, $z \in M$, $Y := (C(M), \leq_C)$. Then $f$ is upper resp. lower semi-continuous at $x \in X$ if and only if for each $\e >0$ and each $R>0$ there is $U \in N(x)$ with $f(u) \cap B (z,R) \subset B^{d_B}(f(x), \e)$ resp. $f(u) \cap B (z,R) \supset B^{d_B}(f(x), \e)$ for all $ u \in U$. 
		\end{enumerate}
\end{Theorem}

\V{\bf Proof.} In Item 1, for lower semicontinuity, choose distinct $x_1, x_2 \in f(x)$, $U_i \in N(x)$ with $f(u) \supset S \setminus \{ x_i \} \forall u \in U_i $. Thus for $U := U_1 \cap U_2 \in N(x) $, we get $f(u) \supset f(x) \forall u \in U$. Similar for each assertion of Items 1 and 2. Item 3 follows from the definition of $d_B$, from $B(X, - \e) := X \setminus (B (X \setminus f(x), \e)) \subset {\rm int} (f(x)) \forall \e >0$, and conversely, $\forall y \in {\rm int } (f(x)) \exists \e >0 : B(X, - \e) \supset y$. \hfill \qed

\bigskip

Another interesting application of semicontinuity is the following: Let $f: X \rightarrow Y$ be a continuous map between two metric spaces. Then 

	$$ \forall p \in Y \forall \e >0 \forall C \subset X {\rm \ compact \ } \exists U \in N(x) \forall q \in U: f^{-1} (q) \cap C \subset B^{g_x}(f^{-1} (p), \e).  $$

Indeed, assume otherwise, then one can construct a sequence $n \mapsto r_n \in C \setminus B(f^{-1}(p))$ with $f(r_n) \rightarrow_{n \rightarrow \infty} p$. But $r$ has a convergent subsequence whose limit $r_\infty \in C \setminus B(f^{-1} (p))$ is mapped by $f$ to $p$, contradiction. 

\bigskip

\V To reformulate this as a semi-continuity property, for a metric space $(Z,d)$ with base point $x_0$, define the {\bf Busemann metric} $d_1$ on the space $C(Z)$ of closed subsets of $Z$ by 

\bean
d_1 = d_1^{x_0, d_0} &:& C(Z) \times C(Z) \rightarrow \R \cup \{ \infty \} , {\rm such \ that, \ for \ all \ } A,B \in C(Z):\\
d_1(A,B) &:&= \sup \{  \vert d_0 ( \{ x \},A) - d_0 (\{ x \},B) \vert \cdot \exp (-d_0( x_0 , x)) ; x \in Z \}  
\eean 

\V Convergence w.r.t. $d_1$ is equivalent to Hausdorff convergence of intersections with every fixed compact subset of $Z$, and if $d$ is Heine-Borel, then $(C(X), d_1)$ is a compact metric space (for details see e.g. \cite{oM-FC}, Sec. 4). It is easy to see that $d_B \vert_{{\rm Gr} (V)} = d_P $ and also $d_B (K,L) = d_{GH} (K \cap \partial B(0,1), L \cap \partial B(0,1))$ for a Euclidean vector space $V$ and two linear subspaces $K,L$ of $V$. So we can rephrase the result above as follows:

\[ \forall p \in Y \forall \e >0 \exists U \in N(x) \exists  A \in B^{d_B} (f^{-1} (p), \e): \forall q \in U : f^{-1} (q) \subset A . \]

\begin{Theorem}
Let $X$ be a Heine-Borel metric space and let $Y$ be a Hausdorff topological space. Then $f^{-1} : Y \rightarrow (C(X), \leq_C)$ is upper semi-continuous.	
\end{Theorem}

\V{\bf Proof.} Via the remark above on $d_1$-convergence, it suffices to show that the inclusion holds w.r.t. the Hausdorff metric after intersecting both sides with compact subsets $C$ of $X$, which has been established above. \hfill \qed

\bigskip

\V {\bf Remark}. In particular we can apply this to the exponential map: Let $(M,g)$ be a Riemannian manifold and $x \in M$. Let $d_g$ be the metric on $M$ induced by $g$, let the Busemann metric $d_B$ on $C(T_xM)$ be induced by $g_x$.
Then $\exp_x^{-1}: (M,d_g) \rightarrow (C(T_xM), d_B)$ is upper semi-continuous. Full continuity is not true: The other possible inclusion fails e.g. in the example of a unit sphere in which $x $ and $p:= -x$ are points opposite to each other: Let $x,p\neq q_n \rightarrow_{n \rightarrow \infty} p$, then $\exp_x^{-1} (q_n) \cap \cl (B(0, \frac{3}{2} \pi))$ all consist of one or two points but $\exp_x^{-1} (p) \cap B(0, \frac{3}{2} \pi)= \partial B(0, \pi)$.

\bigskip

\V Furthermore, we will need some converse of the classical group-theoretical statement of Montgomery-Zippin (\cite{MZ}), accounted for below in Th. \ref{Mont}. A naive guess would be that its hypothesis implies also the converse of its conclusion, i.e. that for any Lie group $G$ and any compact Lie subgroup $K$ there is a neighborhood $U$ of $K$ s.t. for every Lie subgroup $H$ of $G$ contained in $U$ there is $g_U \in G$ with $g_U^{-1}K g_U  $ being a subgroup of $H$ --- but this is obviously wrong, see e.g. the example of $K$ being the compact Abelian group $\R^2 = G = K$, which contains the irrational slope dense subgroups isomorphic to $\R$. Second, one would like to convert the hypothesis, such that $K$ is supposed to be in some neighborhood of the groups $H$, so we need some uniformity. We use the metric and ask: Is it true that for a Lie group $G$ with a bi-invariant Riemannian metric and a compact subgroup $K$ there is $\e >0$ s.t. for every subgroup $H$ of $G$ with $K \subset B(H, \e)$ there is a $g \in G$ with $g^{-1} K g \subset H$? Again, the example of finite subgroups $H$ of $K:= \mathbb{S}^1$ disprove this. The assumption of connectedness of $H$ alone does not help, cf. the connected rational subgroups $H_n$ of $\SSS^1 \times \SSS^1 $ of slope $1/n$ approximating arbitrarily well the entire group. But in this example, every $G$-ball around $1$ contains elements of $H_n$ far from $1$ in the intrinsic distance on $H_n$. This motivates us to include the {\em intrinsic} distance in the hypothesis: 

\begin{Theorem}[Converse Montgomery-Zippin]
\label{MZ-inversion}
Let $G$ be a compact Lie group and let $K$ be a Lie subgroup of $G$. For all $L \in (0; \infty)$, there is $\e \in (0; \infty)$ such that for each Lie subgroup $H$ of $G$ with intrinsic diameter $\diam_H (H) < L $ and

$$  K \subset B_G (H, \e) $$

\V there is $g_{\e} \in B_G(\1, \e)$ with $g_{\e}^{-1} K g_{\e} \subset H $.
\end{Theorem}

\V Preparing a proof of Theorem \ref{MZ-inversion}, the following two theorems (probably already obtained elsewhere) slightly generalize results of Montgomery-Zippin.

\V As in \cite{Kl}, a subset $A$ of a metric space $Z$ is {\bf strongly convex} iff  

\begin{itemize}
	\item For every $p,q \in A$ there is exactly one shortest curve $c_{pq}: [0;1] \mapsto Z$ from $x$ to $y$, and the image of $c_{pq}$ is in $A$.
	\item Every ball in $A$ is convex, i.e., $\forall p,q,x \in A : c_{pq} ([0;1]) \subset B(x, \max \{d(x,p), d(x,q) \})$.
\end{itemize}

\V For each $p \in Z$ we define the {\bf convexity radius} $r_Z(p) $ {\bf at $p$} by $r_Z(p) := \sup \{ \rho >0 \vert B(p, \rho) \ {\rm strongly \ convex}\}$. Then by definition, $r$ is a $1$-Lipschitz function on $Z$. If $Z$ is a Riemannian manifold, then Whitehead's classical result ensures that $r_Z$ is a positive function on $Z$. We define the {\bf convexity radius} $r (Z) := \inf \{ r_Z(x)  \vert x \in Z\}$ {\bf of $Z$}. If $Z$ is a compact Riemannian manifold, a consequence of the facts above is $r (Z) >0$.

\begin{Theorem}
\label{GeodesicInterpolation}
Let $(Z,g) $ be a Riemannian manifold, let $x \in Z$. Then for each two geodesic curves $\g, \d: [0;1] \rightarrow B(x, r(x)/4)$ and for all $a,b \in [0;1]$,

$$d( \gamma(a), \delta (b)) < \max \{ d(\gamma (0) , \delta (0) ), d(\gamma (0), \delta (1)) , d (\gamma (1), \delta (0)), d (\gamma (1), \delta (1))\} .$$

\end{Theorem}

\V {\bf Proof.} We first consider $y \in B(x, r/4)$ and prove $d(y, c(s)) < d(y, c(0)) , d(y, c(1)))$ for all geodesics $c $ with endpoints in $B(x,r/4)$ due to convexity of the ball $B(y,  \max\{ d(y, c(0)), d(y, c(1)) \})$. We apply this to $y$ being the inner points on the geodesics $\gamma$ and $\delta$ and $c$ being the other geodesic in each case. \hfill \qed

\begin{Theorem}
\label{fix}
Let $(Z,g)$ be a Riemannian manifold and let $H \subset {\rm Isom} (M,g) $ be compact. Let $x \in Z $ and let $r$ be the convexity radius of $(Z,g)$ at $x$. If $ H(x) \subset U:= B(x,r/2)$ then $H $ has a fix point in $U$.  
\end{Theorem}

\V {\bf Proof.} Let $A:= \{ p \in Z \vert H(p) \subset \cl (U)\}$. As $x \in A$, the latter is nonempty. It is a closed subset of $\cl (U)$ and thus compact. Let $q \in \cl (U)$ with 

\bea
\label{Minimality}
\diam (H(q)) = \min \{ \diam (H(p)) \vert p \in A \} ,
\eea

\V which exists as $p \mapsto \diam H(p) $ is continuous. Assume that $H(q) \neq \{ q \}$. Then let $\gamma$ be the unique geodesic between two distinct elements $x,z$ of $H(q)$, i.e., $\g: [0;1] \rightarrow U$ with $ \g (0) = x$ and $\g (1) = z$. Choose $a \in (0;1)$ and let $y:= \gamma (a)$. Convexity implies that $y \in A$. We claim that 

\bea
\label{Midpoint}
\diam (H(y)) < \diam (H(q)) ,
\eea 

\V in contradiction with Eq. \ref{Minimality}. And in fact, convexity of $\cl(U)$ implies that $H(y) \subset \cl (U)$ and compactness of $H$ ensures that there is an $h \in H$ with $\diam (H(y)) = d(y, h(y))$. If we define $\d := h \o \g$ (so $\d (a) = h(y))$, then the claim follows by application of Theorem \ref{GeodesicInterpolation} with $a=b$, as $\g (0), \g (1), \d (0) , \d (1) \in H \cdot \g (0) $, so all terms in the maximum are $\leq \diam (H \cdot \g (0))$. \hfill \qed 

\bigskip

\V As a corollary, we recover a version of the Montgomery-Zippin Theorem (cf. \cite{MZ}) with explicit size of the neighborhood in the hypothesis, something we will later need in the proof of the converse version of the same theorem:

\begin{Theorem}[Quantified Montgomery-Zippin Theorem]
\label{Mont}
Let $G$ be a Lie group and let $K$ be a compact Lie subgroup of $G$. Let $r_K$ be the convexity radius of $G/K$ with a left-invariant metric $g$. Let $0< r < r_K$. Then for each Lie subgroup $H$ of $G$ with $H  \subset W:= B(K,r) $, there is $g_{r,H} \in B_G(\1, r)$ s.t. $g_{r,H}^{-1} H g_{r,H} \subset K $.
\end{Theorem}

\V {\bf Proof.} Let $q: G \rightarrow G/K$. First note that $[\1] = \1 \cdot K = K$. Let $H$ be a subgroup of $G$ contained in $W$. Then for all $h \in H$, we have $hK \subset W$, so in the quotient $Z:= G/K$ we have $hK \in U:= \{ w K \vert w \in W  \} \subset G/K$. Thus Theorem \ref{fix} applied to $x = q(\1)$ implies that $H$ has a fix point $gK$ in $U$, i.e. $HgK = gK$, so $g^{-1} Hg K = K$, which implies $ g^{-1} H g \subset K$. \hfill \qed 

\bigskip

\V Now we want to switch roles of $H$ and $K$ in the last step. To this aim, we first restrict ourselves to the case of $G$ carrying a bi-invariant metric, which exists for  $G=\Hol^0$, and more general for each closed, hence compact, subgroup of $SO(n)$, where a bi-invariant metric on $G$ is induced by the negative of the Killing form. Theorem \ref{MZ-inversion} is implied by the following theorem:

\begin{Theorem}
\label{DiamImpliesConv}
Let $g$ be a bi-invariant metric on a Lie group $G$, $L \in (0; \infty)$. Then there is $\e >0$ s.t., for each subgroup $H \neq G$ of $G$ with $\diam (H^0,H^0) < L$, the convexity radius $r_H = r(G/H)$ of $G/H$ satisfies $r_H \geq \e$.
\end{Theorem}

\V {\bf Remark.} The intrinsic diameter bound cannot be omitted, as rational subgroups of tori show. 

\V {\bf Proof.} As the metric is bi-invariant, we can perform for each connected component independently; we focus on the connected component of the identity. The statement will follow e.g. from Dibble's \cite{Dibble} improvement 

$$ r(Z) = \min \{ r_f (Z), \la(Z)/4 \}$$

\V of Klingenberg's Lemma, where $Z$ is a complete Riemannian manifold, $ r_f(Z)$ is the focal radius, which is $\geq \frac{\pi}{\sqrt{\d}}$ if $\sec_Z \leq \d $, and $\la(Z)$ is the length of the shortest nonconstant periodic geodesic in $Z$. The map $q$ as above is a Riemannian submersion. To verify the hypotheses of Dibble's theorem, we first examine $\sec_{G/H}$, which for every, possibly nonconnected, $H$ only depends on the connected component of the identity of $H$. It is well-known (see e.g. \cite{Niko}, Th 3.1) that for a bi-invariant metric $\langle \cdot , \cdot \rangle $, under the identification of horizontal vector fields (i.e. vector fields with values in the orthogonal complement $P:= (T_1H)^\perp$ of $T_1H$) with vector fields in $G/H$ and with $p$ being the orthogonal projection $T_1 G \twoheadrightarrow P$, we get for the Riemannian curvature of $G/H$ and for $X,Y$ horizontal:

\bean
\langle R(X,Y)X, Y \rangle = - \frac{3}{4} \vert p ([X,Y]) \vert^2 - \frac{1}{2} \langle [X, [X,Y]], Y \rangle - \frac{1}{2} \langle [Y, [Y,X]], X \rangle \\
+ \vert U(X,Y)  \vert ^2 - \langle U(X,X) , U(Y,Y) \rangle   
\eean

\V where $U: P \times P \rightarrow P$ is defined by $2 \langle U(X,Y) , Z \rangle  = \langle [Z,X], Y \rangle + \langle X , [Z, Y] \rangle $ for every horizontal $Z$. As $p$ does not increase the length, compactness of the unit sphere implies that there is a global bound on $\sec$.

\bigskip

\V Now for the estimate on $\la(G/H)$: Let $G$ be a fixed $n$-dimensional Lie group. First we show a generalization of well-known results accounted for e.g. in \cite{hB}, Sec. 1.4.

\begin{Lemma}
	\label{OrthogonalFail}
	Let $ a: \N\rightarrow LA(G) \setminus \{ 0\} $ with $a(n) \rightarrow_{n \rightarrow \infty} 0$ and $\exp (a(n)) \in H_n \rightarrow_{n \rightarrow \infty}^{d_{{\rm Hd}}} H_\infty$, and assume that $H_n$ for all $n$ and $H_{\infty}$ are closed subgroups. Let $ a(n) / \vert a (n) \vert \rightarrow_{n \rightarrow \infty} v \in LA(G)$. Then $\exp (t v) \in \exp (V_\infty)$ for all $t \in \R$.  
\end{Lemma}

\V {\bf Proof of the lemma.} Let $t \in \R$, let $ c(n):= \lfloor t/ \vert a(n ) \vert \rfloor \in \N$. Then $ t \in I_n:= [  c(n) \vert a(n) \vert ; (c(n) +1) \vert a(n) \vert )$. As $a$ converges to $0$, the length of the interval $I_n$ converges to zero, thus $c(n) \vert a(n) \vert \rightarrow_{n \rightarrow \infty} t$ and, for all $n \in \N$, 

$$ H_n \ni  (\exp a(n))^{c(n)} =\exp (c(n) \cdot a(n) ) = \exp ( c(n) \vert a(n) \vert \cdot \frac{a(n) }{ \vert a(n) \vert}) \rightarrow_{n \rightarrow \infty} \exp (tv) .$$

\V The assertion of the theorem follows from the above together with the fact $H_\infty = \{ x \in G \vert \exists b : \N\rightarrow G: b(n) \in H_n \forall n \in \N \land b(n) \rightarrow_{n \rightarrow \infty} x\}$ \hfill \qed

\begin{Lemma}
	\label{Closedness}
	Let $D \in [0;\infty]$, $k \in [1;n]$ and let $L_D$ be the space of all Lie subalgebras $V$ of $LA(G)$ such that $\exp (V)$ is of intrinsic diameter $\leq D$. Then $L_D$ is closed in ${\rm Gr}_k(LA(G))$.
\end{Lemma}

\V {\bf Proof of the lemma.} Closedness of the Lie algebra condition is obvious, for the diameter condition let w.l.o.g. $D < \infty$, then let a sequence of Lie subalgebras $n \mapsto V_n \in L_D$ be given that converges to a Lie subalgebra $V_\infty$. Let $\Delta: {Gr}_k (LA(G)) \rightarrow [0;\infty) $ be defined by $ \forall V \in {\rm Gr}_k(LA(G)): \Delta (V) := d_{{\rm Hd}} (\exp (\overline{B}(0,D+1) \cap V), \exp (\overline{B}(0,D) \cap V) ) $. Then $\Delta $ is continuous, and $\Delta (V_n) = 0$ for all $n \in \N$, thus $\Delta (V_\infty) =0$ and thus $\diam (V_\infty) \leq D$. \hfill \qed

\begin{Lemma}
	\label{epsilon}
	 For all $L>0$, there is $\e >0$ such that for all connected Lie subgroups $H$ with intrinsic diameter $< L$ we get 
	 \bea
	 \label{windschief}
	 H \cap B(1, \e)  \subset \exp (B_{T_1 H} (0, \e)).
	 \eea

\end{Lemma}

\V {\bf Proof of the lemma.} The task can be reduced to showing

$$ \exp (B_{T_1 H^\perp }(0,\e))  \setminus \{0 \} ) \cap H = \emptyset. $$

\V Indeed, let $p $ minimize $d(1, \cdot)$ on the compact set $ H \setminus \exp (B_{T_1 H}(0,\e)) \cap \cl (B(1,\e))  $ and assume $p \in \exp (B(0,\e) \setminus (T_1H)^\perp ) $, then the first variational formula for geodesics implies that there is $q \in \exp (T_1 H \cap B(0, \e)) $ with $d(q,p) < d(1,p) < \e$, but $ q^{-1} p \in H$ and $d(1, q^{-1} p ) = d(q,p)$, thus $q^{-1}p \in B(1, \e)$ contradicting the minimizing assumption.

\V Now, assume the opposite of the last displayed equation. Then there is a sequence $a$ of connected Lie subgroups $H$ with intrinsic diameter $< L$ and $w (n) \in T_1 a(n) ^\perp \setminus \{ 0 \}$ such that $\exp (w(n)) \in a(n) $ and $w(n) \rightarrow_{n \rightarrow \infty} 0$. Compactness of $\partial B(0,1)$ implies that the sequence of the $x(n) := w(n) / \vert w(n) \vert$ has a subsequence converging to some $v \in T_1 G$ with $ \vert v \vert = 1 $. Furthermore, by the pidgeon hole principle we can restrict to a subsequence of $k$-dimensional subgroups, for some $k$. Then by compactness of ${\rm Gr}_k (LA(G))$ there is a further subsequence converging to some limit subspace $V_\infty$, which by closedness of the subset of Lie subalgebras is a Lie subalgebra. Furthermore, we have $v \perp V_\infty$ as $v \perp LA(H_n)$ for all $n \in \N$. And $\exp $ maps $V_\infty$ to a closed subgroup $H_\infty$ due to the uniform intrinsic diameter bound on the $V_n$ by Lemma \ref{Closedness}. Then the application of the Lemma \ref{OrthogonalFail} yields $\exp (t v) \in H_\infty \forall t \in \R $, thus $v \in LA (H_\infty)$. But $v \perp LA (H_\infty)$, so $v = 0$ which contradicts $\vert v \vert = 1 $. \hfill \qed

\bigskip

\V By homogeneity we can assume that the geodesic $c$ appearing in $\la(G/H)$ starts at $[\1]$. It lifts to a geodesic $\overline{c}$ of the same length starting at $1$ with differential in $P : = (T_1 H)^\perp$ and with an endpoint in $H$. We can uniformly bound $\{ \vert v \vert : v \perp T_1 H , \exp (v) \in H \}$ from below by the $\e $ from Lemma \ref{epsilon}, which finally concludes the proof of Th. \ref{DiamImpliesConv}. \hfill \qed

\bigskip

\V As now by Th. \ref{DiamImpliesConv} the convexity radius of $G/H$ can be bound uniformly from $0$ for all subgroups $H$ of intrinsic diameter $\leq L$, Th. \ref{MZ-inversion} follows by swapping $K$ and $H$ in the proof of Th. \ref{Mont}. \hfill \qed

\bigskip 

\V To finish the proof of Theorem \ref{main}, let us try to show that the hypotheses of Theorem \ref{MZ-inversion} are satisfied if we approach a limit metric. It turns out that it is essential to take into account the minimal length of curves realizing a given element of the holonomy. Indeed, let $\e >0$. For each $A \in \Hol (g_1)$ there is a curve $c^A$ with $P^{g_1}_{c^A} = A$. Then there is a $C^1$ neighborhood $U$ of $g_1$ with

$$\forall g \in U: A \in B(P_{c^A}^g, \e) \subset B(\Hol(g), \e) .$$ 

\V But the size of the neighborhood depends on the $g_1$-length $l(c_A)$ of $c_A$. 

\V Thus, we consider $\l: \Hol^0(g_1) \rightarrow \R$, $\lambda (A) := \inf \{ l(c) \vert P_c^{g_1} = A\}$. As $\Hol^0(g)$ is compact, a positive lower bound on $\l$ could be established if we knew that $\l$ was upper semi-continuous. Unfortunately, this is wrong i.g., which can be seen by examples of manifolds that are flat in a neighborhood of $p$. (Conversely, we can prove that $\l$ is {\em lower} semi-continuous by an easy application of the general fact that for each $f: A \rightarrow \R$ continuous and $g: A \rightarrow C$ continuous and $f$ bounded on preimages of $g$, we have $\mu: c \mapsto \inf \{ f(a) \vert g(a) = c \} $ is lower semi-continuous). 

\bigskip

\V But there is Wilkins' result \cite{Wilkins} that, whenever $\Hol(g)$ is compact, we get

$$ \forall g \in \Met (M) \exists L_0 \in \R\forall A \in \Hol (g) \exists c \in \Omega (M) : l(c) \leq L_0 \land P^g (c) = A. $$

\V Wilkins' proof can be transferred verbatim to $\Hol^0$ in case of $\Hol^0$ compact (just restrict to contractible loops). Moreover, Wilkins' proof verbatim implies a stronger statement:

\bean \forall g \in \Met (M) \exists L_0 \in \R\forall A \in \Hol (g) \exists c \in \Omega (M) : 
l(c) ,  d^{\Omega_p (M) } (c,p) < L_0 \land P^g (c) = A. 
\eean

  Thus we get via usual ODE estimates that, whenever ${\rm } Hol (g)$ is compact,

\bea
\label{Hol-behavior}
\forall g \in \Met (M) \forall \e >0 \exists \d >0 \forall h \in B_{\Met(M)} (g, \d) : \Hol (g) \subset B_{SO(n)}(B_{\Hol (h)}(\1, \underbrace{L_0 + \e}_{=: L}), \e),  
\eea

\V where $d^{\Omega (M)} (c,p) $ is the distance between $c$ and the constant curve at $p$ in the length space $\Omega(M)$ of loops at $p$ with the intrinsified supremum metric.

\bigskip

\V {\bf Proof of the Main Theorem.} If the diameter of $ \Hol^0 (\nabla ') $ is bounded by $L$ for each $\nabla ' $ in a neighborhood of $\nabla$, by Eq. \ref{Hol-behavior} we know that $ {\rm Hol}^0 (g) \subset B(\Hol^0 (h) , \e)$ and we conclude as desired via Theorem \ref{MZ-inversion} (the converse Montgomery-Zippin Theorem). This shows the first item of the main theorem. For the second item, we use: 

\begin{Theorem}
	\label{UniformDiameterBound}
	$\{ \diam_{\Hol^0 (g)} Hol^0 (g) | M \ n{\rm \ -manifold}, g \in \Riem (M) \}  $ is bounded.	
\end{Theorem}	

\V{\bf Proof.} Go through Berger's finite list and Cartan's classification of symmetric spaces, and noting that there are, up to conjugation, finitely many splittings into invariant subspaces, and that the intrinsic diameter is invariant under conjugation. \hfill \qed

\bigskip

\V {\bf Proof of te Main Theorem.} Taking into account the first item of the Main Theorem shown above, Theorem \ref{UniformDiameterBound} ensures the hypothesis in the case considered in the second item. This concludes the proof of the main theorem. \hfill \qed

\bigskip

\V The choice of $C^2$ regularity of the metric in our main result has been made in view of Wilkin's result which we use and which in turn uses the Ambrose-Singer Theorem. The latter, to our best knowledge, is available if the connection is at least $C^1$ such that curvature exists as a continuous quantity. It is, however, conceivable that the second item of the result is valid even for metrics of regularity $C^{1,1}$ (where we still have local existence and uniqueness of the parallel transport), via the holonomy principle. But not much seems to be known about holonomy in general for metrics of regularity below $C^2$.

\end{document}